\documentclass[a4paper,11pt]{amsart}
\usepackage{amssymb,amsfonts,amsmath}

\def\d{\delta}

\def\C{\mathbb{C}}
\def\c2{\mathbb{C}^2}
\def\R{\mathbb{R}}

\def\1{\bold{1}}
\def\B{\mathbb{B}}

\def\a{\alpha}
\def\b{\beta}
\def\e{\varepsilon}

\def\f{\varphi}
\def\g{\gamma}

\def\p{\psi}

\def\o{\omega}

\newtheorem{lem}{Lemma}[section]
\newtheorem{pro}[lem]{Proposition}

\newtheorem{def/not}[lem]{Definition/Notations}

\newtheorem{thm}[lem]{Theorem}

\newtheorem{exa}[lem]{Example}

\begin{document}

\title[H\"older continuous solutions to Monge-Amp\`ere equations]
{H\"older continuous solutions to Monge-Amp\`ere equations}

\author{V.GUEDJ \& S.KOLODZIEJ \& A.ZERIAHI}

\begin{abstract}
We study the regularity of solutions to complex Monge-Amp\`ere equations
$(dd^c u)^n=f dV$,
on bounded strongly pseudoconvex domains $ \Omega \subset \C^n$.
We show, under a mild technical assumption, that the unique solution
$u$ to such an equation is H\"older continuous if the boundary
values $\phi$ are H\"older continuous and the 
density $f$ belongs to $L^p(\Omega)$ for some $p>1$.
This improves previous results by Bedford-Taylor and Kolodziej.
\end{abstract}

\maketitle


\section*{Introduction}

Let $\Omega$ be a bounded strongly pseudoconvex open subset of $\C^n$.
Given $\phi \in {\mathcal C}^0(\partial \Omega)$ and $f \in L^p(\Omega)$,
we consider the Dirichlet problem
$$
MA(\Omega,\phi,f): \left\{ \begin{array}{l}
 u \in PSH(\Omega) \cap {\mathcal C}^0(\overline{\Omega}) \\
 u=\phi \text{ on } \partial \Omega \\
 (dd^c u)^n=f \beta_n \text{ in } \Omega.
 \end{array} \right.
$$
Here $\b_n=dV$ denotes the euclidean volume form in $\C^n$, $d=\partial+\overline{\partial}$,
$d^c=i(\overline{\partial}-\partial)$, $PSH(\Omega)$ is the set of plurisubharmonic 
(psh for short) functions
in $\Omega$ (the set of locally integrable functions $u$ such that $dd^c u \geq 0$
in the sense of currents), and $( dd^c \; \cdot )^n$ denotes the complex Monge-Amp\`ere operator:
this operator is well defined on the subset of bounded (in particular continuous) psh functions,
as follows from the work of E.Bedford and A.Taylor [BT 2]. We refer the reader to
[K 4] for a recent survey on its properties.

The equation $MA(\Omega,\phi,f)$ has been studied intensively during the last decades.
Let ${\mathcal B}(\Omega,\phi,f)$ denote the family of subsolutions to 
$MA(\Omega,\phi,f)$, i.e. the set of bounded functions $v$ that are psh in $\Omega$ with
$v \leq \phi$ on $\partial \Omega$ (i.e. $\limsup_{\zeta \rightarrow z} v(\zeta) \leq \phi(z)$
when $z \in \partial \Omega$) and $(dd^c v)^n \geq f \b_n$
in $\Omega$. It follows from the comparison principle [BT 2] that if a solution exists,
it must coincide with the Perron-Bremermann envelope,
$$
u(z)=u(\Omega,\phi,f)(z):=\sup \{ v(z) \, / \, v \in {\mathcal B}(\Omega,\phi,f) \}.
$$

It follows from the works of H.J.Bremermann [Br], J.Walsh [W] 
and E.Bedford-A.Taylor [BT 1] 
that $u(\Omega,\phi,f)$ is indeed
a solution to $MA(\Omega,\phi,f)$ when $f \in {\mathcal C}^0(\overline{\Omega})$.
S.Kolodziej has further shown in [K 2] 
that the solution $u(\Omega,\phi,f)$
is still continuous when $f \in L^p(\Omega)$, for some
$p>1$ (it is easy to check that $u(\Omega,\phi,f)$ is not necessarily
locally bounded when the density is merely in $L^1(\Omega)$).

Higher regularity results have also been provided. E.Bedford and A.Taylor showed in
[BT 1] that if $\phi \in Lip_{2\a}(\Omega)$ and $f^{1/n} \in Lip_{\a}(\overline{\Omega})$,
then $u(\Omega,\phi,f) \in Lip_{\a}(\overline{\Omega})$.
The smoothness of the solution (assuming smoothness of $\phi$ and $f>0$) is established in [CKNS].
Our aim here is to establish H\"older-continuity of $u(\Omega,\phi,f)$ in the spirit
of [BT 1], only assuming the density $f$ belongs to $L^p(\Omega)$, $p >1$,
as in [K 2]. Our main result is the following:

\vskip.2cm \noindent 
{\bf Theorem A.}
{\it Assume $f \in L^p(\Omega)$, for some $p>1$, and $\phi \in Lip_{2\a}(\partial \Omega)$,
with $\nabla u(\Omega,\phi,0) \in L^2(\Omega)$. Then
$$
u(\Omega,\phi,f) \in Lip_{\a'}(\overline{\Omega}), 
\; \; 
\text{ for all } \a'<\min(\a,2/[qn+2]), 
$$
where $1/p+1/q=1$.}
\vskip.2cm

The condition $\nabla u(\Omega,\phi,0) \in L^2(\Omega)$ is automatically
satisfied if $\phi \in {\mathcal C}^{1,1}(\partial \Omega)$: in this
case $u(\Omega,\phi,0) \in Lip_1(\overline{\Omega})$, hence
$\nabla u(\Omega,\phi,0)$ is actually bounded in $\Omega$ (see [BT 1]).
What really matters here is that there should exist a subsolution 
$v \in {\mathcal B}(\Omega,\phi,0)$ such that $\nabla v \in L^2(\Omega)$.
This implies (see Lemma 3.1), that 
$u(\Omega,\phi,0)$ and $u(\Omega,\phi,f)$
both have gradient in $L^2(\Omega)$.

We could not avoid the use of this extra (technical ?) hypothesis on
the homogenous solution $u(\Omega,\phi,0)$.
Also the exponent $\a'$ is probably not optimal. We can get a better
exponent by assuming that $\Delta u(\Omega,\phi,0)$ has finite mass in $\Omega$
(this is automatically satisfied when $\phi \in {\mathcal C}^2 (\partial \Omega)$).

\vskip.2cm \noindent 
{\bf Theorem B.}
{\it Assume $f \in L^p(\Omega)$, for some $p>1$, and 
$\phi \in Lip_{2\a}(\partial \Omega)$ is such that
that $\Delta  u(\Omega,\phi,0)$ has finite mass in $\Omega$. 
Then
$$
u(\Omega,\phi,f) \in Lip_{\a''}(\overline{\Omega}), 
\; \; 
\text{ for all } \a''<\min(\a,2/[qn+1]), 
$$
where $1/p+1/q=1$.}
\vskip.2cm

The exponent $\a''$ is not far from being optimal as Example 4.2 shows.

\section{The stability estimate}

Our main tool is the following estimate which is proved in [EGZ]
in a compact setting (under growth -- but no boundary -- conditions, see
Proposition 3.3 in [EGZ]).
A similar -- but weaker -- estimate was established by S.Kolodziej in [K 3].

\begin{thm}
Fix $0 \leq f \in L^p(\Omega)$, $p>1$.
Let $\f,\p$ be two bounded psh functions in $\Omega$ such that
$
(dd^c \f)^n=f \beta_n \text{ in } \Omega, 
\text{ and } \f \geq \p \text{ on } \partial \Omega.
$
Fix $r \geq 1$ and $0 \leq \g <r/[nq+r]$, $1/p+1/q=1$. Then 
$$
\sup_{\Omega} (\p-\f) \leq C ||(\max(\p-\f,0)||_{L^r(\Omega)}^{\g},
$$
for some uniform constant $C=C(\g,||f||_{L^p(\Omega)},||\p||_{L^{\infty}})>0$.
\end{thm}

The proof follows closely the one given in [EGZ], so we only
recall the main ingredients, for the reader's convenience.
The estimate is a simple consequence of the following result.

\begin{pro}
Fix $f \in L^p(\Omega)$, $p>1$, and let
$\f,\p$ be bounded psh functions in $\Omega$ such that
$\f \geq \p$ on $\partial \Omega$.
If $(dd^c \f)^n=f \beta_n$, then for all $\e>0,\tau>0$,
$$
\sup_{\Omega} (\p-\f) \leq \e+C \left[ \text{Cap}(\f-\p<-\e ) \right]^{\tau},
$$
for some uniform constant $C=C(\tau, ||f||_{L^p(\Omega)})$.
\end{pro}

Here $\text{Cap}(\cdot)$ denotes the Monge-Amp\`ere capacity introduced and studied
by E.Bedford and A.Taylor in [BT 2].
Recall that for $K \subset \Omega$,
$$
\text{Cap}(K):=\sup \left\{ \int_K (dd^c v)^n \, / \, v \in PSH(\Omega)
\text{ with } -1 \leq v \leq 0 \right\}.
$$

The proposition is a direct consequence of the following three lemmas:

\begin{lem}
Fix $\f,\p \in PSH(\Omega) \cap L^{\infty}(\Omega)$ such that 
$\overline{\lim}_{\zeta \rightarrow \partial \Omega} (\f-\p) \leq 0$.
Then for all $t,s>0$, 
$$
t^n \text{Cap}(\f-\p<-s-t) \leq \int_{(\f-\p<-s)} (dd^c \f)^n.
$$
\end{lem}

\begin{lem}
Assume $0 \leq f \in L^p(\Omega)$, $p>1$.
Then for all $\tau>1$, there exists $C_{\tau}>0$ such that
for all $K \subset \Omega$,
$$
0 \leq \int_K f dV \leq C_{\tau} \left[ \text{Cap}(K) \right]^{\tau}.
$$
\end{lem}

When $(dd^c \f)^n=f dV$, one can 
combine Lemma 1.3 and Lemma 1.4 to control
$\text{Cap}(\f-\p<-s)$ by $\text{Cap}(\f-\p<-s-t)]^{\tau}$.
This has strong consequences since $\tau>1$, as the following
result shows, when applied to $g(t)=\text{Cap}(\f-\p<-t-\e)^{1/n}$:

\begin{lem}
Let $g:\R^+ \rightarrow \R^+$ be a decreasing right-continuous function.
Assume there exists $\tau,B>1$ such that
$g$ satisfies
$$
H(\a,B) \hskip1cm
tg(s+t) \leq B [ g(s) ]^{\tau}, \;
\forall s,t>0.
$$

Then there exists $S_{\infty}=S_{\infty}(\tau,B) \in \R^+$ such that
$g(s)=0$ for all $s \geq S_{\infty}$.
\end{lem}

One word about the proofs. Lemma 1.3 is a direct consequence
of the ``comparison principle'' of E.Bedford and A.Taylor [BT 2]
(see [K 4] p.32 for a detailed proof).
Using H\"older's inequality, one reduces the proof of Lemma 1.4 to 
showing that the euclidean volume is bounded from above by the 
Monge-Amp\`ere capacity. One can actually show that
$$
Vol(K) \lesssim \exp[-Cap(K)^{-1/n}] \; \;  
(\text{see Theorem 7.1 in [Z]}),
$$
which is a much better control than what we actually need.
The last lemma is an elementary exercise, whose proof is given
in [EGZ], Lemma 2.3.

\section{H\"older continuous barriers}

For fixed $\d >0$ we consider
$\Omega_{\d}:=\{ z \in \Omega \, / \, \text{dist}(z,\partial \Omega) >\d \}$ and set
$$
u_{\d}(z):=\sup_{||\zeta|| \leq \d } u(z+\zeta), \, z \in \Omega_{\d}.
$$
This is a psh function in $\Omega_{\d}$, when $u$ is psh in $\Omega$,
 which measures the modulus of continuity
of $u$. We would like to use Theorem 1.1 applied with $\p=u_{\d}$.
However $u_{\d}$ is not globally defined in $\Omega$, so we need
to extend it with control on the boundary values.
This is the contents of our next result which makes heavy use of the
pseudoconvexity assumption.

\begin{pro} Let $u \in PSH (\Omega) \cap L^{\infty} (\Omega)$ be a psh 
function such that $u_{|\partial \Omega}=\phi \in Lip_{2 \a}(\partial \Omega)$.
Then there exists a family $(\tilde{u}_{\d})_{0 < \d < \d_0}$ 
of bounded psh functions on $\Omega$ 
such that $\tilde{u}_{\d} \searrow u$ in $\Omega$ as $\d \searrow 0$, with
$$
\tilde{u}_{\d} = \left\{ \begin{array}{l}
\max (u_{\d} - C \d^{\a}, u)  \text{ in } \Omega_{\d} \\
u  \text{ in } \Omega \setminus \Omega_{\d}
\end{array} \right.
$$
\noindent In particular $\sup_{\Omega_{\d}} \vert \tilde{u}_{\d} - u_{\d}\vert \leq C \d^{\a}$ 
for $0 < \d < \d_0$.
\end{pro} 

The proof relies on the construction of H\"older continuous plurisubharmonic 
and plurisuperhamonic barriers for the Dirichlet problem $MA (\Omega,\phi,f)$:

\begin{lem} Fix $\phi \in Lip_{2 \a} (\partial \Omega)$, $f \in L^p(\Omega)$,
$p>1$, and set $u := u(\Omega,\phi,f)$. 
Then there exists $v, w \in PSH (\Omega) \cap Lip_{\a} (\overline \Omega)$ such that 
\begin{enumerate}
 \item $ v (\zeta) = \phi (\zeta) = - w (\zeta), \forall \zeta \in \partial \Omega$,\\

 \item $ v (z) \leq u (z) \leq - w (z), \forall z \in \Omega$.
\end{enumerate}
\end{lem}

\begin{proof}
Assume first that $\phi \equiv 0$. We are going to show that
there exists a weak barrier $b_f \in PSH (\Omega) \cap Lip_1 (\Omega)$ for the Dirichlet problem 
$MA (0,f,\Omega)$, i.e. a psh function which satisfies
\begin{itemize}
\item $(i)$ $ b_f (\zeta) = 0, \  \forall \zeta \in \partial \Omega,$\\

\item $(ii)$ $ b_f \leq u (\Omega,0,f), $ in $\Omega$,\\

\item $(iii)$ $  \vert b_f (z) - b_f (\zeta) \vert \leq C_1 \vert z - \zeta \vert, \ \forall z \in \Omega, \ \forall \zeta \in \Omega$, 
\end{itemize}
for some uniform constant  $C_1 > 0$.

In order to construct $b_f$, we set $u_0 := u (\Omega,0,f)$ and assume first that 
the density $f$ is bounded near $\partial \Omega$: there 
exists a compact subset $K \subset \Omega$ such that $0 \leq f \leq M$ on $\Omega \setminus K$.
Let $\rho$ be a ${\mathcal C}^2$ strictly plurisubharmonic defining function for $\Omega$.
Then for $A > 0$ large enough the function $b_f := A \rho$ satisfies the condition 
$(dd^c b_f)^n \geq M \beta_n \geq f \beta_n$ on $\Omega \setminus K$. 
Moreover taking $A$ large enough we also have $ A  \rho \leq  m \leq u_0$ on a neighborhood of $K$, 
where $m := \min_{\Omega} u_0.$
Therefore the function $b_f$ is a ${\mathcal C}^2$ plurisubharmonic function on $\Omega$ satisfying the 
conditions $(dd^c b_f)^n \geq (dd^c u_0)^n$ on $\Omega \setminus K$ and 
$b_f \leq u_0$ on $\partial (\Omega \setminus K)$. 
This implies, by the comparison principle [BT 2], 
that $b_f \leq u_0$ in $\Omega \setminus K$, hence in $\Omega$.

When $f$ is not bounded near $\partial \Omega$, we can proceed as follows.
Fix a large ball $\B \subset \C^n$ so that $\Omega \Subset \B \subset \C^n$. 
Define $\tilde f := f$ in $\Omega$ and $\tilde f = 0$ in $\B \setminus \Omega$. 
We can use our previous construction to find a barrier function 
$b_{\tilde f} \in PSH (\B) \cap {\mathcal C}^2 (\B)$ for the Dirichlet problem 
$MA (\B,0, \tilde f)$ for the ball $\B$. 
Let  $h= u (\Omega,-b_{\tilde f},0)$ denote the Bremermann function
in $\Omega$ with boundary values $-b_{\tilde f}$, for the zero density.
Since $-b_{\tilde f} \in {\mathcal C}^2(\partial \Omega)$,
the psh function $h$ is Lipschitz in $\Omega$ (see [BT 1]),
therefore $b_f :=h + b_{\tilde f} \in PSH (\Omega) \cap Lip_1 (\Omega)$ 
is a barrier function for the Dirichlet problem $MA(\Omega,0,f)$. 

It remains to construct the functions $v,w$ satisfying (1),(2) above.
It follows from [BT 1] that the psh
functions $u (\Omega,\pm \phi,0)$ are H\"older continuous of order $\a$. 
We let the reader check that the functions $v := u (\Omega,\phi,0) + b_f$ 
and $w := u (\Omega,-\phi,0) + b_f$ do the job.
\end{proof}

We are now ready for the proof of the proposition.

\begin{proof}
It follows from Lemma 2.2 that
$$
|u(z) - u (\zeta)| \leq C |z - \zeta|^{\a},  
\; \forall \zeta \in \partial \Omega, \forall z \in \Omega.
$$
The functions $u_{\d}(z):=\sup_{||\zeta|| \leq \d} u(z+\zeta)$   
are psh in  $\Omega_{\d}$.
Observe that if $z \in \partial \Omega_{\delta}$ and $\zeta  \in \C^n$ with 
$||\zeta||\leq \delta$ then  $z + \zeta \in \partial \Omega$, hence
$u_{\delta} - C \delta^{\a} \leq u (z)$. 
Thus the functions 
$$ 
\tilde{u}_{\delta} (z) := \left\{ \begin{array}{l}
\sup \{u_{\delta} (z) - C \delta^{\a}, u (z)\} \text{ in } \Omega_{\delta} \\
u \text{ in } \Omega \setminus \Omega_{\delta} \end{array} \right.
$$
are psh and bounded in  $\Omega$
and decrease to $u$ as $\d$ decreases to $0$.
\end{proof}

Our construction of barriers allows us to control the total
mass of the Laplacian of solutions to $MA(\Omega,\phi,f)$.
This will be important in section 4.

\begin{pro}
Fix $0 \leq f \in L^p(\Omega)$, $p>1$, and $\phi \in Lip_{2 \a}(\partial \Omega)$.

Then $\Delta u(\Omega,0,f)$ has finite
mass in $\Omega$.
In particular, if $\Delta u(\Omega,\phi,0)$ has finite mass in $\Omega$,
then $\Delta u(\Omega,\phi,f)$ has finite mass in $\Omega$.
\end{pro}

Note that $\Delta u(\Omega,\phi,0)$ has finite mass in $\Omega$ when
$\phi \in {\mathcal C}^2(\Omega)$,
as explained in the proof below.

\begin{proof}
Assume first that $\phi \in {\mathcal C}^2(\Omega)$.
Consider any smooth extension of $\phi$ and correct it by adding $A \rho$, $A>>1$,
in order to obtain a smooth plurisubharmonic extension $\hat{\phi}$
which is defined in
a neighborhood of $\overline{\Omega}$. 
Since $\hat{\phi}$ is a subsolution to $MA(\Omega,\phi,0)$ whose
Laplacian has finite mass in $\Omega$, it follows from the comparison
principle that $\Delta u(\Omega,\phi,0)$ also has finite mass
in $\Omega$.

Let $\tilde{f}$ be the trivial extension of $f$ to a large ball $\B$
containing $\Omega$. Let $b_{\tilde{f}} \in {\mathcal C}^2(\B)$ 
be a plurisubharmonic barrier for $MA(\B,0,\tilde{f})$ (see the proof of Lemma 2.2).
Then $b_f:=u(\Omega,-b_{\tilde{f}},0)+b_{\tilde{f}}$ is
a plurisubharmonic barrier for $MA(\Omega,0,f)$. Its Laplacian has finite mass
in $\Omega$ since $b_{\tilde{f}}$ is smooth, so it follows from the comparison
principle that $\Delta u(\Omega,0,f)$ has finite mass in $\Omega$.

Set now $v:=u(\Omega,0,f)+u(\Omega,\phi,0)$.
This is a plurisubharmonic function in $\Omega$ such that
$v=\phi$ on $\partial \Omega$ and $(dd^c v)^n \geq f dV$ in $\Omega$.
If $\Delta u(\Omega,\phi,0)$ has finite mass in $\Omega$, then
$\Delta v$ has finite mass in $\Omega$, hence
$\Delta u(\Omega,\phi,f)$ also has finite mass in $\Omega$.
\end{proof}

\section{Gradient estimates}

This section is devoted to the proof of Theorem A.
Let $u=u(\Omega,\phi,f)$ be the unique solution to the 
complex Monge-Amp\`ere equation
$$
(dd^c u)^n=f \beta_n \text{ in } \Omega,
$$
with boundary values $u=\phi \in Lip_{2\a}(\partial \Omega)$. 
Since $f \in L^p(\Omega)$, $p>1$, it follows from [K 2] that 
$u$ is a continuous plurisubharmonic function.
Our aim is to show that $u$ is H\"older continuous.

Let $\tilde{u}_{\d}$ be the functions given by Proposition 3.1.
The stability estimate (Theorem 1.1) applied with $r=2$ yields
$$
\sup_{\Omega_{\d}} (u_{\d}-u) \leq
C_1 \d^{\a}+\sup_{\Omega} (\tilde{u}_{\d}-u) 
\leq C_1 \d^{\a}+ C_2 ||u_{\d}-u||_{L^2(\Omega_{\d})}^{\g},
$$
for $\g<2/(nq+2)$, $1/p+1/q=1$.
It remains to
show that $||u_{\d}-u||_{L^2(\Omega_{\d})} \leq C_4 \d$ to conclude the proof.

It will be a consequence of Lemma 3.1 below that $\nabla u \in L^2(\Omega)$.
Assuming this for the moment, we derive the desired upper-bound
on $||u_{\d}-u||_{L^2(\Omega_{\d})}$.
Averaging the gradient of u on balls of radius $\d$ yields, for all $||\zeta||<\delta$,
$$
\left( \int_{\Omega_{\d}} |u(z+\zeta)-u(z)|^2 dV(z) \right)^{1/2}
\leq \d \; || \nabla u||_{L^2(\Omega)}.
$$
By Choquet's lemma, there exists a sequence $\zeta_j$, $||\zeta_j|| <\d$, such that
$u_{\d}=( \sup_j u_j )^*$, in $\Omega_{\d}$,
where $u_j(z):=u(z+\zeta_j)$.
Since $u_j-u \leq u_{\d}-u$, 
it follows from Lebesgue's dominated convergence theorem that
$$
\left( \int_{\Omega_{\d}} |u_{\d}(z)-u(z)|^2 dV(z) \right)^{1/2}
\leq \d \; || \nabla u||_{L^2(\Omega)}.
$$
This ends the proof of Theorem A up to the fact, to be established now, that
$u$ has gradient in $L^2(\Omega)$.

\vskip.2cm

Since $u$ is plurisubharmonic and continuous, $\nabla u \in L_{loc}^2(\Omega)$.
It follows from Lemma 3.1 below that $\nabla u \in L^2(\Omega)$ as soon as
$u$ is bounded from below by a continuous psh function $v$ such that
$v \leq u \text{ in } \Omega$, $v=u=\phi$ on $\partial \Omega$, and $\nabla v \in L^2(\Omega)$.
Our extra assumption in Theorem A precisely yields such a function $v$.
Indeed set $v:=u(\Omega,\phi,0)+b_f$, where $b_f$ is the psh barrier 
constructed in the proof of Lemma 2.2: this is a psh function such that
\begin{itemize}
\item $v=\phi+0=u$ on $\partial \Omega$;
\item $(dd^c v)^n \geq (dd^c b_f)^n \geq f \beta_n$ in $\Omega$, thus $v \leq u$ in $\Omega$;
\item $\nabla u(\Omega,\phi,0) \in L^2(\Omega)$ 
and $\nabla b_f \in L^{\infty}(\Omega)$, hence $\nabla v \in L^2(\Omega)$.
\end{itemize}

It is easy to check that $\nabla u(\Omega,\phi,0) \in L^{\infty}(\Omega) \subset L^2(\Omega)$
when $\phi \in {\mathcal C}^2(\partial \Omega)$.
We refer the reader to  [BT 1] for a proof of the more delicate result that this still holds
when $\phi \in {\mathcal C}^{1,1}(\partial \Omega)$.

\begin{lem} Let $u, v \in PSH (\Omega) \cap C^{0} (\overline{\Omega})$ such that
$v \leq u$ on $ \Omega$ and $v = u$ on $\partial \Omega$. Then 
$ \int_{\Omega} d u \wedge d^c u \wedge \omega^{n - 1} 
\leq  \int_{\Omega} d v \wedge d^c v \wedge \omega^{n - 1}.$
\end{lem}

\begin{proof}
First assume  that $u = v$ near the boundary $\partial \Omega$. 
Then we can approximate by smooth plurisubharmonic functions $u_{\delta}\downarrow u$ and 
$v_{\delta}  \downarrow v$ with $u_{\delta} \leq v_{\delta}$ on $\Omega_{\delta}$ and 
$u_{\delta} = v_{\delta}$ 
in a neighborhood of $\partial  \Omega_{\e},$ for $0 < \delta < \e$ and $\e$ small enough. 
Integrating by parts, we get
$$
\int_{\Omega_{\e}} d v_{\delta} \wedge d^c v_{\delta} \wedge \omega^{n - 1} =  
\int_{\partial \Omega_{\e}} v_{\delta} d^c v_{\delta} \wedge \omega^{n - 1} -
\int_{\Omega_{\e}} v_{\delta} dd^c v_{\delta} \wedge \omega^{n - 1} .
$$
Since $v_{\delta} = u_{\delta}$ on a neighbourhood of $\partial \Omega_{\e}$ we conclude that 
$v_{\delta} d^c v_{\delta} = u_{\delta} d^c u_{\delta}$ on this neighbourhood and 
then using again integration by parts we get
$$
\int_{\Omega_{\e}} d v_{\delta} \wedge d^c v_{\delta} \wedge \omega^{n - 1} =
\int_{\Omega_{\e}}  d u_{\delta} \wedge d^c u_{\delta} \wedge \omega^{n - 1}
+ \int_{\Omega_{\e}} (u_{\delta} dd^c u_{\delta} - v_{\delta} dd^c v_{\delta}) \wedge \omega^{n - 1}. 
$$
On the other hand, since $v_{\delta} = u_{\delta}$ on a neighbourhood of 
$\partial \Omega_{\e}$ we conclude that
$v_{\delta} \wedge d^c u_{\delta} = u_{\delta} \wedge d^c v_{\delta}$ on this neighbourhood and then 
\begin{eqnarray*}
\int_{\Omega_{\e}} u_{\delta} dd^c v_{\delta} - v_{\delta} dd^c u_{\delta}) \wedge \omega^{n - 1}
& = &\int_{\partial \Omega_{\e}} (u_{\delta} d^c v_{\delta} - v_{\delta} d^c u_{\delta}) \wedge \omega^{n - 1} \\
& - & \int_{\Omega_{\e}} (d u_{\delta} \wedge d^c v_{\delta} - d v_{\delta} 
\wedge d^c u_{\delta}) \wedge \omega^{n - 1} = 0.
\end{eqnarray*}
Therefore
$$
\int_{\Omega_{\e}} d v_{\delta} \wedge d^c v_{\delta} \wedge \omega^{n - 1}
=  \int_{\partial \Omega_{\e}} d u_{\delta} \wedge d^c u_{\delta} \wedge \omega^{n - 1}
 + \int_{\Omega_{\e}} (u_{\delta} - v_{\delta}) (dd^c u_{\delta} + dd^c v_{\delta}) \wedge \omega^{n - 1}.
$$
 Since $v_{\delta} \leq u_{\delta}$ on $\Omega_{\delta},$ we obtain 
$$ 
\int_{\Omega_{\e}} d v_{\delta} \wedge d^c v_{\delta} \wedge \omega^{n - 1}
\geq \int_{\partial \Omega_{\e}} d u_{\delta} \wedge d^c u_{\delta} \wedge \omega^{n - 1}.
$$
We know by Bedford and Taylor's convergence theorem that
$d v_{\delta} \wedge d^c v_{\delta} \wedge \omega^{n - 1} 
\rightarrow d v \wedge d^c v \wedge \omega^{n - 1}.$
Taking the limit when $\delta \searrow 0$ we get
$$
\int_{\overline{\Omega}_{\e}} d v \wedge d^c v \wedge \omega^{n - 1}
\geq \int_{\Omega_{\e}} d u \wedge d^c u \wedge \omega^{n - 1}.
$$
Taking the limit when $\e \downarrow 0$ we get the required inequality.

 Now if we only know that $u = v,$ we can define for $\e > 0$
$u_{\e} := \sup \{u - \e , v\}.$ Then $v \leq u_{\e}$ on $\Omega$ and $u_{\e} = v$ 
near the boundary of $\Omega$. Therefore we have for $\delta > 0$ small enough 
$$
\int_{\Omega} d v \wedge d^c v \wedge \omega^{n - 1} \geq \int_{\Omega} d u_{\e} 
\wedge d^c v_{\e} \wedge \omega^{n - 1}.
$$
 Now by Bedford and Taylor's convergence theorem, we know that
$d u_{\e} \wedge d^c u_{\e} \wedge \omega^{n - 1} \to d u \wedge d^c u \wedge \omega^{n - 1}$ 
as $\e \downarrow 0$. Therefore we have
$$
\int_{\Omega} d v \wedge d^c v \wedge \omega^{n - 1} \geq \int_{\Omega} 
d u \wedge d^c u \wedge \omega^{n - 1},
$$
which proves the required inequality. 
\end{proof}

\section{Laplacian estimates}

This section is devoted to the proof of Theorem B.
We use the same method as above. The finiteness of the
total mass of $\Delta u(\Omega,\phi,0)$
allows a good control (see Lemma 4.2)
on the terms $\hat u_{\delta} - u$,
where
$$
 \hat u_{\delta} (z) := \frac{1}{\tau_{n} \delta^{2 n}} \int_{\vert 
 \zeta - z\vert \leq \delta} u (\zeta) d V (\zeta), z \in \Omega_{\delta},
$$
where $\tau_n$ denotes the volume of the unit ball in $\C^n$.
We shall compare $\hat u_{\d}$ to $u_{\d}$ in Lemma 4.1 below.

It follows from the construction of plurisubharmonic H\"older continuous barriers 
that the solution $u = u (\Omega,\phi,f)$ is H\"older continuous near the boundary,
i.e. for $\delta > 0$ small enough, we have
\begin{equation}
u (z) - u (\zeta) \leq c_0 \delta^{\a},
\end{equation}
for $z, \zeta \in \overline{\Omega}$ with 
$\hbox{dist} (z, \partial \Omega) \leq \delta, \hbox{dist} (\zeta, \partial \Omega) \leq \delta $ 
and $ \vert z - \zeta \vert \leq \delta.$
\vskip.1cm

The link between $u_{\d}$ and $\hat{u}_{\d}$ is made by the
following lemma.

\begin{lem} Given $\alpha \in ]0, 1[$, the following two conditions are equivalent. 

$(i)$ There exists $\delta_0, A > 0$ such that for any 
$0 < \delta \leq \delta_0$,
$$ 
u_{\delta} - u \leq A \delta^{\alpha} \, \hbox{ on } \, \Omega_{\delta}.
$$

$(ii)$ There exists $\delta_1, B > 0$ such that for any $0 <\delta< \delta_1$,
$$
\hat u_{\delta} - u \leq  B \delta^{\alpha} \, \hbox{ on }  \Omega_{\delta}.
$$
\end{lem}

\begin{proof}
Observe  that $\hat u_{\delta} \leq u_{\delta}$ in $\Omega_{\delta}$, hence
$(i) \Longrightarrow (ii)$ follows immediately.

We now prove that $(ii) \Longrightarrow (i)$.  
We need to show that there exists $A,\d_0>0$ such that
for $0 < \delta \leq \delta_0$,
$$
\omega (\delta) := \sup_{z \in \Omega_{\delta}} [u_{\delta} (z) - u (z)] \leq A \delta^{\alpha}.
$$
Fix $\d_{\Omega}>0$ small enough so that $\Omega_{\d} \neq \emptyset$ for 
$\d \leq 3 \d_{\Omega}$.
Since $u$ is uniformly continuous, for any fixed
$0< \delta <\d_{\Omega}$,
$$ 
\nu (\delta) := \sup_{\delta < t \leq \delta_{\Omega}}   
\omega (t)  t^{- \alpha} < + \infty. 
$$

We claim that there exists $\delta_0  > 0$  small enough so that for any  
$0 < \delta  \leq \delta_0$,
$$ 
\omega (\delta) \leq  A \delta^{\alpha}, 
\text{ with }
A= (1+4^{\a}) c_0+2^{\a} 4^n B+\nu(\d_{\Omega}),
$$
where $c_0$ is the constant arising in inequality (1), while
$B$ is the constant from condition (ii).
Assume this is not the case. 
Then there exists $0 < \delta  < \delta_{\Omega}$ such that
\begin{equation}
\omega (\delta)> A \delta^{\alpha}.
\end{equation}
Set $\d:=\sup \{ t <\d_{\Omega} \, / \, \o(t)> A t^{\a} \}$.
Then 
\begin{equation}
\frac{\o(\d)}{\d^{\a}} \geq A \geq \frac{\o(t)}{t^{\a}}
\text{ for all } t \in [\d,\d_{\Omega}].
\end{equation}
Since $u$ is continuous, we can find
$z_0 \in \overline{\Omega_{\delta}}$, $  \zeta_0 \in \overline{\Omega}$ with 
$ \vert z_0 - \zeta_0 \vert  \leq  \delta$ s.t.
$$
\omega(\delta) = \sup_{z \in \Omega_{\d}} \left[ \sup_{w \in B(z,\d)} u(w)  -u(z) 
\right]=u (\zeta_0) - u (z_0).
$$
\vskip.2cm

We first derive a contradiction if $z_0$
is close enough to the boundary of $\Omega$.
Assume that $\hbox{dist} (z_0,\partial \Omega) \leq 3  \delta$. 
Take $z_1 \in \partial \Omega$ such that
$\hbox{dist} (z_0,\partial \Omega)  = \hbox{dist} (z_0, z_1) \leq 4   \delta$.
It follows from $(1)$ that
$$ 
\omega(\delta) =  u (\zeta_0) - u (z_0) = 
[u (\zeta_0) - u (z_1)] + [u (z_1) - u (z_0)] \leq [1+4^{\a}] c_0 \delta^{\a}.
$$
This contradicts $(3)$.

Thus we can assume that $\hbox{dist} (z_0,\partial \Omega) >  3  \delta$.
Fix $b>1$ so that $\hbox{dist} (z_0,\partial \Omega) >  (2b+1)  \delta$.
Thus any $z \in \B(\zeta_0,b\d)$ satisfies 
$z \in \B(z_0,[b+1]\d)$, hence $z \in \Omega_{b \d}$.
By using inequality (3) with $t=b\d$, we get
$u(\zeta_0)-u(z) \leq b^{\a} \o(\d)$, hence
\begin{equation}
u(z) \geq u(\zeta_0)-b^{\a} \o(\d), 
\, \text{ for all } z \in \B(\zeta_0,b\d).
\end{equation}

Observe now that $\B(\zeta_0,\d) \subset \B(z_0,[b+1]\d)$, hence
\begin{eqnarray*}
\hat u_{(b+1) \delta} (z_0) 
& = & \left(\frac{b}{b+1} \right)^{2n} \hat u_{b \delta} (\zeta_0)  +   
\frac{1}{\tau_n (b+1)^{2n} \delta^{2 n} }
\int_{ \B (z_0, (b+1) \delta) \setminus  \B (\zeta_0, b \delta)} u d V \\
& \geq &  
\left(\frac{b}{b+1} \right)^{2n}  u(\zeta_0)  + 
\left[(1 - \frac{b^{2n}}{(b+1)^{2 n}}\right] 
[u (\zeta_0) -  b ^{\alpha} \omega (\delta)] \\
&=& u(\zeta_0) -b^{\a} \left[1 - \frac{b^{2n}}{(b+1)^{2 n}}\right] \o(\d),
\end{eqnarray*}
where we have used the subharmonicity of $u$ together with inequality (4).
Since $u(\zeta_0)=u(z_0)+\o(\d)$, we infer, letting $b \rightarrow 1$,
$$
\hat u_{2 \delta} (z_0)  \geq u(z_0)
+4^{-n}\o(\d).
$$
 
We now use assumption (ii), only considering small enough values of $\d>0$:
since $\hat{u}_{2 \d}(z_0) \leq u(z_0)+B 2^{\a} \d^{\a}$, we get
$$
\o(\d) \leq 4^n 2^{\a} B \d^{\a} < A \d^{\a}.
$$
This contradicts the definition of $\d$, hence we have proved that 
$(ii) \Rightarrow (i)$.
\end{proof}

It is straightforward to check that if (i) is satisfied, then
$u$ belongs to $Lip_{\a}(\overline{\Omega})$.
Thus Theorem B will be proved if we can establish (ii).
It follows from Theorem 1.1 that is suffices to get control
on the $L^1$-average of $\hat{u}_{\d}-u$.
This is the contents of our next result.

\begin{lem}
For $\delta > 0$ small enough, we have
$$
\int_{\Omega_{\delta}} [\hat u_{\delta} (z) - u (z)] d V_{2 n} (z) \leq c_n 
\Vert \Delta u\Vert \delta^2, 
$$
where $c_n > 0$ is a uniform constant.
\end{lem}
 
\begin{proof}
It follows from Jensen's formula that for $z \in \Omega_{\delta}$ and $0 < r < \delta$,
$$ 
\frac{1}{\sigma_{2 n - 1}} \int_{\vert \xi\vert = 1} u (z + r \xi) d S_{2 n - 1} = u (z) + 
\int_0^{r} t^{1 - 2 n} (\int_{\vert \zeta\vert \leq t} dd^c u \wedge \beta_{n - 1}) d t.
$$
Using polar coordinates we get, for $z \in \Omega_{\delta}$,
$$
 \hat u_{\delta} (z) - u (z)  = \frac{1}{\sigma_{2 n - 1} 
\delta^{2 n}} \int_{0}^{\delta} r^{2 n - 1} d r \int_0^{r} t^{1 - 2 n} 
(\int_{\vert \zeta - z\vert \leq t} dd^c u \wedge \beta_{n - 1}) d t.
$$
Finally Fubini's theorem yields
\begin{eqnarray*}
\int_{\Omega_{\delta}} (\hat u_{\delta} - u) d V_{2 n} &\leq&
 a_n \delta^{- 2 n}\int_{0}^{\delta} r^{2 n - 1} d r 
\int_0^{r} t^{1 - 2 n} (\int_{\vert \zeta\vert \leq t} (\int_{\Omega} \Delta u) d t \\
&\leq& c_n  \delta^2 \Vert \Delta u\Vert.
\end{eqnarray*}
\end{proof}

This ends the proof of Theorem B since by proposition 2.3, 
$\Delta u=\Delta u(\Omega,\phi,f)$ has finite mass in $\Omega$.
\vskip.2cm

We now give a simple example which shows that one can not expect a better 
exponent than $\a=2/nq$, for $1/p+1/q=1$.

\begin{exa}
Consider $u(z_1,\ldots,z_n):=|z_1|^{\a} \cdot ||(z_2,\ldots,z_n)||^2$.
This is a plurisubharmonic function in $\C^n$ which is H\"older-continuous
of exponent $\a \in ]0,1[$. We let the reader check that
$$
(dd^c u)^n=f dV, \; 
\text{ with } f(z)=\frac{1}{|z_1|^{2-n\a}} g(z_2,\ldots z_n),
$$
where $g>0$ is a smooth density.

Given $p>1$, $f$ belongs to  $L_{loc}^p(\C^n)$ whenever
$\a=\e+2/nq$, for some $\e>0$.
This shows that we cannot get a better exponent than 
$2/nq$ in Theorems A,B.
\end{exa}

\vskip .2cm

Slawomir Kolodziej

Jagiellonian University

Institute of Mathematics

Reymonta 4, 30-059 KRAKOW (POLAND)

Slawomir.Kolodziej@im.uj.edu.pl

\vskip.2cm

Vincent Guedj \& Ahmed Zeriahi

Laboratoire Emile Picard

UMR 5580, Universit\'e Paul Sabatier

118 route de Narbonne

31062 TOULOUSE Cedex 04 (FRANCE)

guedj@picard.ups-tlse.fr ; zeriahi@picard.ups-tlse.fr

\end{document}